\documentclass[11pt]{article}
\usepackage{graphicx,amssymb,mathrsfs,amsmath}
\newfont{\bb}{msbm10}

\begin{document}
\cleardoublepage \pagestyle{myheadings}

\bibliographystyle{plain}

\title{\bf On the generalized low rank approximation of the correlation matrices arising in the asset portfolio
          \thanks{\scriptsize The work was supported by National Natural Science Foundation of China (Nos. 11101100; 11261014;
11301107; 61362021), Natural Science Foundation of Guangxi Province (No. 2012GXNSFBA053006;
2013GXNSFBA019009; 2013GXNSFBB053005; 2013GXNSFDA019030), the Fund for Guangxi Experiment
Center  of  Information  Science  (20130103),  Innovation  Project  of  GUET  Graduate  Education
(GDYCSZ201473), Innovation Project of Guangxi Graduate Education (YCSZ2014137), and Guangxi Key
Lab of Wireless Wideband Communication and Signal Processing open grant 2012.} \\}
\author
{Xuefeng Duan \footnote{\scriptsize Corresponding author.E-mail address:duanxuefenghd@aliyun.com(X. Duan), baijianchaok@126.com(J. Bai).}
\quad Jianchao Bai
\quad Maojun Zhang
\quad Xinjun Zhang  \\
{\small\it  College of Mathematics and Computational Science,
Guilin University of}\\
{\small\it Electronic Technology, Guilin 541004, P.R. China} \\}
\date{}
\maketitle

\hrule
\bigskip
{\noindent\bf Abstract} \vskip 1mm
\small
In this paper, we consider the generalized low rank approximation of the correlation matrices problem
which arises in the asset portfolio. We first characterize the feasible set by using
the Gramian representation together with a special trigonometric function transform, and then transform
the generalized low rank approximation of the correlation matrices problem into an unconstrained
optimization problem. Finally, we use the conjugate gradient algorithm with the strong Wolfe line
search to solve the unconstrained optimization problem. Numerical examples show that
our new method is feasible and effective.


\bigskip
\hrule
\bigskip
\noindent{\normalsize\bf 1. Introduction} \vskip2mm

Throughout this paper, we use $R^{n\times n}$ and $S^+_n$ to denote the set of $n\times n$ real matrices
and symmetric positive semidefinite matrices, respectively. We use $A^T$ and $tr(A)$ to represent the
transpose and trace of the matrix $A$, respectively. The symbols $\|A\|_F$ and $rank(A)$ denote the Frobenius norm
and the rank of the matrix $A,$ respectively. The symbol $diag(Y)$ stands for the vector whose elements lie
in the diagonal line of the matrix  $Y,$ and the symbol $e$ stands for the vector whose elements are of all
ones, i.e., $e=(1,1,\cdots,1)^{T}.$

In this paper, we consider the following problem named {\bf generalized low rank
approximation of the correlation matrices}.

\noindent{\bf Problem 1.1.} Given some correlation matrices $A^{(d)}\in R^{n\times n},\ d=1,2,\cdots,m$,
and a positive integer $k,\ 1\leq k <n$, find a correlation matrix $\widehat{Y}$ whose rank is less than
and equal to $k$ such that $$\frac{1}{2}\sum\limits_{d=1}^{m}\|A^{(d)}-\widehat{Y}\|^2_F
=\min\limits_{Y\in S^+_n,diag(Y)=e,rank(Y)\leq k}\frac{1}{2}\sum\limits_{d=1}^{m}\|A^{(d)}-Y\|^2_F. \eqno(1.1)$$

Problem (1.1) arises in the asset portfolio (see [10] for more details), which can be stated as follows.
Suppose that $R=DCD$ is the covariance matrix of $n$ assets, where $C$ is a correlation matrix
and $D$ is a diagonal matrix with positive variances which are specially used to describe
the risk of assets. In practice, the covariance matrix is usually estimated by the
historical data of the return of each asset, that is, an approximation covariance is obtained
by statistics method. Let $$R^{(d)}=D^{(d)}C^{(d)}D^{(d)}$$
be the approximation covariance with $d$th sampling some data,
where $D^{(d)}$ and $C^{(d)}$ are the $d$th approximation diagonal matrix and correlation matrix, respectively.
Higham [4] proposed a method for finding the nearest low rank approximation of a correlation matrix by only one
sampling(i.e., $m=1$). However, it is difficult for the decision maker to choose the best approximation covariance
matrix with only one sampling because there is always a noise in the data on the prices of assets. Thus, we develop
a repeated sampling method to get a series of approximation covariance matrices, that is, $d$ comes from $1$ to $m$.
Obviously, it is very easy to obtain the optimal diagonal matrix $\widehat{D}$ by a series of $D^{(d)}$.
The major obstacle to finding the optimal covariance matrix is conducting the optimal correlation
matrix $\widehat{C}$ from a series of $C^{(d)}$. The above consideration leads to solving the following problem:
given some correlation matrices $A^{(1)},A^{(2)},\cdots,A^{(m)}\in R^{n\times n}$, find a correlation
matrix $\widehat{Y}$ such that $$\frac{1}{2}\sum\limits_{d=1}^{m}\|A^{(d)}-\widehat{Y}\|^2_F=\min\limits_{Y\in S^+_n,\ diag(Y)=e}\frac{1}{2}\sum\limits_{d=1}^{m}\|A^{(d)}-Y\|^2_F. \eqno(1.2)$$
Meanwhile, for the large financial correlation matrices, usually almost all variances can be attributed to some
stochastic Brownian factors. Therefore, instead of taking into account all Brownian motions, we would wish to
simulate with a smaller number of factors, i.e., $rank(Y) < n$ and typically $rank(Y)$ is from 1 to $k$. Then the
problem (1.2) with rank constraint becomes problem (1.1).

Noting that the matrix $Y$ in problem (1.1)  is not only positive semidefinite but also satisfies $rank(Y)\leq k$,
so problem (1.1) belongs to the structured low rank approximation problem. As Gillard-Zhigljavsky [3] said, the
structured low rank approximation is a difficult optimization problem, so there is much work to be done.
In the last few years, there has been a constantly increasing interest in developing the theory and numerical
methods for the nearest low rank approximation of a correlation matrix, due to their wide applications in the
fiance and risk management [6], machine learning [15], stress testing of bank [13], industrial process
monitoring [7] and image processing [5]. Recently, problem (1.1) with $m=1$ has been extensively studied,
and the research results mainly concentrate on the following two cases. One is without the rank constraint
and the other is with the rank constraint.

For the case without the rank constraint, Higham [4] proposed an alternative projection
algorithm to solve the nearest correlation matrix problem by defining two projection operators.
Under some proper assumptions, Li-Li [8] developed a projected semismooth Newton method to solve the problem
of calibrating least squares covariance matrix. Qi and Sun [12] proposed a Newton-type method for the nearest
correlation matrix problem, and the quadratic convergence of the new method was proved. An unconstrained convex
optimization approach was proposed to find the nearest correlation matrix to the target matrix with the
fixed correlations unaltered in [13].  Besides, Qi-Sun [14] introduced an augmented Lagrangian dual method
for for the H-weighted nearest correlation matrix problem. This method solves a sequence of unconstrained strongly convex optimization
problems, each of which can be solved by a semismooth Newton method combined with the conjugate gradient method.
Recently, Yin, etc [18, 20] developed two new alternative gradient algorithms to compute the nearest correlation
matrix by making use of the alternative gradient method.

For the case with the rank constraint, by making use of the fact that
$$ Y\in S^+_n,\ rank(Y)\leq k \Longleftrightarrow \lambda_{k+1}(Y)+\cdots+\lambda_{n}(Y)=0, $$
Gao and Sun [2] proposed a majorized penalty approach for solving the rank constrained correlation matrix problem.
It is noted that Gao and Sun's  majorized penalty approach can deal with  some large scale problems ($n\geq 500$).
Motivated by the method in [12] and based on a well-known result that the sum of the largest
eigenvalues of a symmetric matrix can be represented as a semidefinite programming problem,
Li-Qi [9] proposed a novel sequential semismooth Newton method to solve problem (1.1) with $m=1$.
They formulate the problem as a bi-affine semidefinite programming and then use an augmented
Lagrange method to solve a sequence of least squares problems.
Both Simon-Abell [16] and Pietersz-Groenen [11] used majorization approach
to solve the low rank approximation of a correlation matrix. The difference lies in that the former solved
the problem with any weighted norm while the latter only settled it with Frobenius norm.
By constructing a Lagrange function, Zhang-Wu [21] transformed the
low rank approximation of a correlation matrix into a min-max problem, where the inner maximization problem
was solved with closed form spectral decomposition and the outer minimization problem was solved with
gradient-based methods. In [1], Grubisic and Pietersz introduced a geometric programming approach
to solve the low rank nearest correlation matrix problem. The method could be used to minimize any sufficiently
smooth objective function.

However, the research results of problem (1.1) with $m>1$ are very few as far as we know. The greatest
difficulties to solve problem (1.1) are how to characterize the feasible set and deal with the complex
structure. In this paper, we overcome these difficulties by using the Gramian representation together
with a special trigonometric function transform. Then problem (1.1) is transformed into an unconstrained
optimization problem. Finally, the conjugate gradient method with the strong Wolfe line search is
given to solve the unconstrained optimization problem. Numerical examples show that our
new method is feasible and effective.

\vskip 3mm \noindent{\normalsize\bf 2. Main results }
\vskip2mm

In this section, we first transform problem (1.1) into an unconstrained optimization problem
by making use of the Gramian representation together with a special trigonometric function transform.
Then we use the conjugate gradient algorithm with the strong Wolfe line search to solve it.

We first define the following set $$S=\{Y\in R^{n\times n} \mid \ Y\in S^+_n,\ rank(Y)\leq k\}.$$
It is easy to characterize the set $S$ by using the Gramian representation (see [17]), i.e., $$Y=XX^T,\ \ X\in R^{n\times k}.$$

Set $$\Gamma=\{Y\in R^{n\times n} \mid \  diag(Y)=e\}.$$  It is easy  to verify that the feasible
set of problem (1.1) is $S\bigcap \Gamma$. The most difficulty to solve problem (1.1) is how to
characterize the feasible set. Now we begin to use the Gramian representation together with a
special trigonometric function transform to characterize the feasible set $S\bigcap \Gamma$.

\noindent{\bf Theorem 2.1.} Let the matrix $X$ be
$$X=[X_1,X_2,\cdots,X_k]=\left[\begin{array}{cccc}
x_{11} & x_{12} & \cdots & x_{1k} \\
x_{21} & x_{22} & \cdots & x_{2k} \\
\vdots & \vdots & \ddots & \vdots \\
x_{n1} & x_{n2} & \cdots & x_{nk}
\end{array}\right]\in R^{n\times k}.$$
Suppose
$$X_1=\left[\begin{array}{c}
\cos\alpha_{11} \\
\cos\alpha_{21} \\
\vdots \\
\cos\alpha_{n1}
\end{array}\right],\
\ X_2=\left[\begin{array}{c}
\cos\alpha_{12}\sin\alpha_{11} \\
\cos\alpha_{22}\sin\alpha_{21} \\
\vdots \\
\cos\alpha_{n2}\sin\alpha_{n1}
\end{array}\right],\cdots,$$

$$X_{k-1}=\left[\begin{array}{c}
\cos\alpha_{1k-1}\prod\limits_{l=1}^{k-2}\sin\alpha_{1l}\\
\cos\alpha_{2k-1}\prod\limits_{l=1}^{k-2}\sin\alpha_{2l} \\
\vdots \\
\cos\alpha_{nk-1}\prod\limits_{l=1}^{k-2}\sin\alpha_{nl}
\end{array}\right],\
\ X_{k}=\left[\begin{array}{c}
\prod\limits_{l=1}^{k-1}\sin\alpha_{1l}\\
\prod\limits_{l=1}^{k-1}\sin\alpha_{2l} \\
\vdots \\
\prod\limits_{l=1}^{k-1}\sin\alpha_{nl}
\end{array}\right],$$
where $\alpha_{ij}\in R,\ i=1,2,\cdots,n,\ j=1,2,\cdots,k-1,$
then the matrix $Y=XX^T\in R^{n\times n}$ is not only symmetric positive semidefinite,
but also satisfies $rank(Y)\leq k$ and $diag(Y)=e.$

\noindent{\bf Proof}. By using the Gramian representation, it is easy to verify that the
matrix $Y$ is symmetric positive semidefinite and satisfies $rank(Y)\leq k$. Hence, we only need to prove $diag(Y)=e$.

Consider the matrix $X$ with $k=2$. According to the assumptions, we have
$$X=[X_1,X_2]=\left[\begin{array}{cc}
\cos\alpha_{11} & \sin\alpha_{11}\\
\cos\alpha_{21} & \sin\alpha_{21}\\
\vdots & \vdots\\
\cos\alpha_{n1}& \sin\alpha_{n1}
\end{array}\right].$$
Let $\chi_{i}\ (i=1,2,\cdots,n)$ be the $i$th row of the matrix $X$, that is,
$$\chi_i=[\cos\alpha_{i1} , \sin\alpha_{i1}].$$ By multiplying $\chi_i$ and $\chi_i^T$,
we get the element $y_{ii}$ of the matrix $Y$, that is,
$$y_{ii}=\chi_i\cdot \chi_i^T=(cos\alpha_{i1})^2+(sin\alpha_{i1})^2=1.$$
That is to say, $diag(Y)=e$. Hence, Theorem 2.1 holds when $k=2$.

When $k>2$, without loss of generality, we take the $i$th row of the matrix $X$ and write
it as $\chi_i, \ i=1,2,\cdots,n$, then
$$\chi_i=[\cos\alpha_{i1},\sin\alpha_{il}\cos\alpha_{i2},\cdots,
\cos\alpha_{ik-1}\prod\limits_{l=1}^{k-2}\sin\alpha_{il},\prod\limits_{l=1}^{k-1}\sin\alpha_{il}].$$
By multiplying $\chi_i$ and $\chi_i^T$, we get the element $y_{ii}$ of the matrix $Y$, that is,
$$ \begin{array}{lll}
y_{ii}&=&\chi_i\cdot \chi_i^T\\ &=&(cos\alpha_{i1})^2+(\sin\alpha_{il}\cos\alpha_{i2})^2
+\cdots+(\cos\alpha_{ik-1}\prod\limits_{l=1}^{k-2}\sin\alpha_{il})^2
+(\prod\limits_{l=1}^{k-1}\sin\alpha_{il})^2  \\
&=&(cos\alpha_{i1})^2+(\sin\alpha_{il}\cos\alpha_{i2})^2+\cdots
+(\prod\limits_{l=1}^{k-2}\sin\alpha_{il})^2(\cos^2\alpha_{ik-1}+\sin^2\alpha_{ik-1})\\
&=&(cos\alpha_{i1})^2+(\sin\alpha_{il}\cos\alpha_{i2})^2+\cdots+(\cos\alpha_{ik-2}\prod\limits_{l=1}^{k-3}\sin\alpha_{il})^2
+(\prod\limits_{l=1}^{k-2}\sin\alpha_{il})^2\\
&=&(cos\alpha_{i1})^2+(\sin\alpha_{il}\cos\alpha_{i2})^2+\cdots
+(\prod\limits_{l=1}^{k-3}\sin\alpha_{il})^2(\cos^2\alpha_{ik-2}+\sin^2\alpha_{ik-2})  \\
&=& (cos\alpha_{i1})^2+(\sin\alpha_{il}\cos\alpha_{i2})^2+\cdots+(\cos\alpha_{ik-3}\prod\limits_{l=1}^{k-4}\sin\alpha_{il})^2
+(\prod\limits_{l=1}^{k-3}\sin\alpha_{il})^2\\
&=& \cdots \\
&=& (cos\alpha_{i1})^2+(\sin\alpha_{il}\cos\alpha_{i2})^2+(\sin\alpha_{il}\sin\alpha_{i2})^2\\
&=&(cos\alpha_{i1})^2+(sin\alpha_{i1})^2\\
&=&1.
\end{array} $$
Hence, for any $k\geq 2$, we have $y_{ii}=1,\ i=1,2,\cdots,n$, that is, $diag(Y)=e.$ $\ \ \ \ \Box$

\noindent{\bf Remark 2.1}. As Simon and Abell [16] said, a correlation matrix is a symmetric
positive semidefinite matrix with unit diagonal, and any symmetric positive semidefinite matrix with unit
diagonal is a correlation matrix. In Theorem 2.1, the matrix $Y$ must be a correlation matrix, and noting that
 $\alpha_{ij},\ \ i=1,2,\cdots, n, \ \ j=1,2,\cdots, k-1$ are arbitrary real number, so the matrix $Y=XX^{T}$
 can be represented all the correlation matrices.

\noindent{\bf Remark 2.2}. To explain Theorem 2.1, we take a $3\times 2$ matrix for example. Set
$$X=[X_1,X_2]=\left[\begin{array}{cc}
cos\alpha_{11} & sin\alpha_{11}  \\
cos\alpha_{21} & sin\alpha_{21}  \\
cos\alpha_{31} & sin\alpha_{31}
\end{array}\right].$$

By a simple calculation, we can obtain that
$$ \begin{array}{lll}
Y&=&XX^T\\
&=&\left[\begin{array}{cc}
cos\alpha_{11} & sin\alpha_{11}  \\
cos\alpha_{21} & sin\alpha_{21}  \\
cos\alpha_{31} & sin\alpha_{31}
\end{array}\right]
\left[\begin{array}{cc}
cos\alpha_{11} & sin\alpha_{11}  \\
cos\alpha_{21} & sin\alpha_{21}  \\
cos\alpha_{31} & sin\alpha_{31}
\end{array}\right]^T \\\\
&=&
\left[\begin{array}{ccc}
1 & cos\alpha_{11}cos\alpha_{21}+sin\alpha_{11}sin\alpha_{21} & cos\alpha_{11}cos\alpha_{31}+sin\alpha_{11}sin\alpha_{31}\\
cos\alpha_{11}cos\alpha_{21}+sin\alpha_{11}sin\alpha_{21} & 1  & cos\alpha_{21}cos\alpha_{31}+sin\alpha_{21}sin\alpha_{31}\\
cos\alpha_{11}cos\alpha_{31}+sin\alpha_{11}sin\alpha_{31} & cos\alpha_{21}cos\alpha_{31}+sin\alpha_{21}sin\alpha_{31} & 1
\end{array}\right].
\end{array} $$
Obviously, the matrix $Y$ is not only symmetric positive semidefinite , but also satisfies $rank(Y)\leq 2$ and $diag(Y)=e.$

By using the similar way in the proof of Theorem 2.1, we can obtain the other elements of the
matrix $Y$, that is,
$$ Y=(y_{ij})_{n\times n}=
\left\{\begin{array}{cc}
\sum\limits_{p=1}^{k-1}cos\alpha_{ip}cos\alpha_{jp}\prod\limits_{l=1}^{p-1}\sin\alpha_{il}\sin\alpha_{jl}
+\prod\limits_{l=1}^{k-1}\sin\alpha_{il}\sin\alpha_{jl},& i\neq j\\
1, & i=j
\end{array}\right..$$

Substituting $y_{ij}$ into problem (1.1), it is easy to obtain that problem (1.1) can be
written as the following unconstrained optimization problem.

\noindent{\bf Problem 2.1}. Given some correlation matrices $A^{(d)}=(A^{(d)}_{ij})_{n\times n},\ d=1,2,\cdots,m$,
and a positive integer $k,\ 1\leq k <n$, find the solution $\widehat{\alpha}\in R^{n\times(k-1)}$ of the following optimization problem
$$\min\limits_{\alpha\in R^{n\times(k-1)}} F(\alpha),\eqno(2.1) $$
where $$F(\alpha)=\sum\limits_{d=1}^{m}\sum\limits_{i=1}^{n-1}\sum\limits_{j=i+1}^{n}(\sum\limits_{p=1}
^{k-1}cos\alpha_{ip}cos\alpha_{jp}\prod\limits_{l=1}^{p-1}\sin\alpha_{il}\sin\alpha_{jl}
+\prod\limits_{l=1}^{k-1}\sin\alpha_{il}\sin\alpha_{jl}-A^{(d)}_{ij})^2. \eqno(2.2)$$

Nextly, we will use the conjugate gradient algorithm with the strong Wolfe line search to solve the
 unconstrained optimization problem. The most difficulty to solve problem (2.1) is how to
compute the gradient of the objective function $F(\alpha)$. Now we begin to compute the
gradient of the objective function.

\noindent{\bf Theorem 2.2}. The gradient of the objective function $F(\alpha)$ of problem (2.1) is
$$\nabla F(\alpha)=(\frac{\partial F(\alpha)}{\partial \alpha_{11}},\frac{\partial F(\alpha)}
{\partial \alpha_{21}},\cdots,\frac{\partial F(\alpha)}{\partial \alpha_{n1}},\cdots,\frac{\partial F(\alpha)}
{\partial \alpha_{1k-1}},\frac{\partial F(\alpha)}{\partial \alpha_{2k-1}},\cdots,\frac{\partial F(\alpha)}{\partial \alpha_{nk-1}})^T,$$
where
$$ \begin{array}{lll} \frac{\partial F(\alpha)}{\partial \alpha_{\mu\nu}}
&=& 2\sum\limits_{d=1}^{m}\sum\limits_{i=1,i\neq \mu}^{n}
\{(\sum\limits_{p=1}^{k-1}cos\alpha_{\mu p}cos\alpha_{ip}\prod\limits_{l=1}^{p-1}\sin\alpha_{\mu l}\sin\alpha_{il}
+\prod\limits_{l=1}^{k-1}\sin\alpha_{\mu l}\sin\alpha_{il}-A^{(d)}_{\mu i})\\
&\times&
(-\sin\alpha_{\mu \nu}\cos\alpha_{i\nu}\prod\limits_{l=1}^{\nu-1}\sin\alpha_{\mu l}\sin\alpha_{i l}+
\cos\alpha_{\mu\nu}sin\alpha_{i\nu}\prod\limits_{l=1,l\neq\nu}^{k-1}\sin\alpha_{\mu l}\sin\alpha_{il}
\\
&+&\cos\alpha_{\mu\nu}\sin\alpha_{i\nu}\sum\limits_{p=\nu+1}^{k-1}cos\alpha_{\mu p}cos\alpha_{ip}
\prod\limits_{l=1,l\neq\nu}^{p-1}\sin\alpha_{\mu l}\sin\alpha_{i l})
\}, \end{array} \eqno(2.3)$$
here $\mu=1,2,\cdots,n, \  \nu=1,2,\cdots,k-1$.

\noindent{\bf Proof}. To prove Theorem 2.2, we only need to prove (2.3) holds when $m=1$,
because the forms of the expression of the gradient of the objective function $F(\alpha)$
with $m=1$ are the same as that with $m>1$.

For $m=1$, noting that the total numbers including $\alpha_{\mu\nu}$ in $F(\alpha)$ are
$$\sum\limits_{i=1}^{\mu-1}(\sum\limits_{p=1}^{k-1}cos\alpha_{\mu p}cos\alpha_{i p}\prod\limits_{l=1}^{p-1}\sin\alpha_{\mu l}\sin\alpha_{il}
+\prod\limits_{l=1}^{k-1}\sin\alpha_{\mu l}\sin\alpha_{il}-A_{i\mu})^2+$$
$$\sum\limits_{j=\mu+1}^{n}(\sum\limits_{p=1}^{k-1}cos\alpha_{\mu p}cos\alpha_{j p}\prod\limits_{l=1}^{p-1}\sin\alpha_{\mu l}\sin\alpha_{jl}
+\prod\limits_{l=1}^{k-1}\sin\alpha_{\mu l}\sin\alpha_{jl}-A_{\mu j})^2. $$
Hence, the derivative of $F(\alpha)$ at $\alpha_{\mu\nu}$ is
$$ \begin{array}{lll} \frac{\partial F(\alpha)}{\partial \alpha_{\mu\nu}}
&=&\frac{\partial}{\partial \alpha_{\mu\nu}} \{  \sum\limits_{i=1}^{\mu-1}(\sum\limits_{p=1}^{k-1}cos\alpha_{\mu p}cos\alpha_{i p}\prod\limits_{l=1}^{p-1}\sin\alpha_{\mu l}\sin\alpha_{il}
+\prod\limits_{l=1}^{k-1}\sin\alpha_{\mu l}\sin\alpha_{il}-A_{i\mu})^2\\
&+& \sum\limits_{j=\mu+1}^{n}(\sum\limits_{p=1}^{k-1}cos\alpha_{\mu p}cos\alpha_{j p}\prod\limits_{l=1}^{p-1}\sin\alpha_{\mu l}\sin\alpha_{jl}
+\prod\limits_{l=1}^{k-1}\sin\alpha_{\mu l}\sin\alpha_{jl}-A_{\mu j})^2\} \\
&=& 2\sum\limits_{i=1}^{\mu-1}\{(\sum\limits_{p=1}^{k-1}cos\alpha_{\mu p}cos\alpha_{i p}\prod\limits_{l=1}^{p-1}\sin\alpha_{\mu l}\sin\alpha_{il}
+\prod\limits_{l=1}^{k-1}\sin\alpha_{\mu l}\sin\alpha_{il}-A_{i\mu})\\
&\times&
(-\sin\alpha_{\mu \nu}\cos\alpha_{i\nu}\prod\limits_{l=1}^{\nu-1}\sin\alpha_{\mu l}\sin\alpha_{i l}+
\prod\limits_{l=1,l\neq\nu}^{k-1}\sin\alpha_{\mu l}\sin\alpha_{il}cos\alpha_{\mu\nu}sin\alpha_{i\nu}
\\
&+&\sum\limits_{p=\nu+1}^{k-1}cos\alpha_{\mu p}cos\alpha_{ip}\prod\limits_{l=1,l\neq\nu}^{p-1}\sin\alpha_{\mu l}\sin\alpha_{i l}\cos\alpha_{\mu\nu}\sin\alpha_{i\nu}
)\}\\
&+& 2\sum\limits_{j=\mu+1}^{n}\{(\sum\limits_{p=1}^{k-1}cos\alpha_{\mu p}cos\alpha_{j p}\prod
\limits_{l=1}^{p-1}\sin\alpha_{\mu l}\sin\alpha_{jl}
+\prod\limits_{l=1}^{k-1}\sin\alpha_{\mu l}\sin\alpha_{jl}-A_{\mu j})\\
&\times&
(-\sin\alpha_{\mu \nu}\cos\alpha_{j\nu}\prod\limits_{l=1}^{\nu-1}\sin\alpha_{\mu l}\sin\alpha_{j l}+
\prod\limits_{l=1,l\neq\nu}^{k-1}\sin\alpha_{\mu l}\sin\alpha_{jl}cos\alpha_{\mu\nu}sin\alpha_{j\nu}
\\
&+&\sum\limits_{p=\nu+1}^{k-1}cos\alpha_{\mu p}cos\alpha_{jp}\prod\limits_{l=1,l\neq\nu}^{p-1}\sin\alpha_{\mu l}\sin\alpha_{j l}\cos\alpha_{\mu\nu}\sin\alpha_{j\nu}
)\}. \end{array} $$
Because $A_{i\mu}=A_{\mu i}$, we turn $j$ to $i$ and conclude that
$$ \begin{array}{lll}
\frac{\partial F(\alpha)}{\partial \alpha_{\mu\nu}}
&=& 2\sum\limits_{i=1,i\neq\mu}^{n}\{(\sum\limits_{p=1}^{k-1}cos\alpha_{\mu p}cos\alpha_{i p}\prod
\limits_{l=1}^{p-1}\sin\alpha_{\mu l}\sin\alpha_{il}
+\prod\limits_{l=1}^{k-1}\sin\alpha_{\mu l}\sin\alpha_{il}-A_{\mu i})\\
&\times&
(-\sin\alpha_{\mu \nu}\cos\alpha_{i\nu}\prod\limits_{l=1}^{\nu-1}\sin\alpha_{\mu l}\sin\alpha_{i l}+
\cos\alpha_{\mu\nu}sin\alpha_{i\nu}\prod\limits_{l=1,l\neq\nu}^{k-1}\sin\alpha_{\mu l}\sin\alpha_{il}
\\
&+&\cos\alpha_{\mu\nu}\sin\alpha_{i\nu}\sum\limits_{p=\nu+1}^{k-1}cos\alpha_{\mu p}cos\alpha_{ip}
\prod\limits_{l=1,l\neq\nu}^{p-1}\sin\alpha_{\mu l}\sin\alpha_{i l}
)\},
\end{array}$$
where $\mu=1,2,\cdots,n,\ \nu=1,2,\cdots,k-1$. $\ \ \ \ \Box$

Consequently, the conjugate gradient algorithm with the strong Wolfe line search to solve the
 minimization problem (2.1) can be described in Algorithm 2.1.

\vskip4mm
\hrule\vskip2mm
\noindent {\bf Algorithm 2.1} (This algorithm attempts to solve problem (2.1))\\
\indent Step 1. Given parameters $\rho \in (0,1),\ \delta\in (0,0.5), \sigma\in(\delta,0.5)$,
and tolerance error $ 0\leq tol\ll 1$. Choose an initial iterative matrix $\alpha_0\in R^{n\times(k-1)}$. Set $t:=0$. \\
\indent Step 2. Calculate $g_t=\nabla F(\alpha_t)$. If $\parallel g_t\parallel_F< tol$, stop
and output $\alpha^*\approx\alpha_t$.\\
\indent Step 3. Determine the search direction $d_t$, where
$$ d_t=
\left\{\begin{array}{cc}
-g_t,& t=0\\
-g_t+\frac{g^T_tg_t}{g^T_{t-1}g_{t-1}}d_{t-1}, & t\geq 1
\end{array}\right.. $$\\
\indent Step 4. Confirm the step length $\beta_t$ by applying the strong Wolfe line search, i.e.,
$$ \left\{\begin{array}{c}
F(\alpha_{t+1})\leq F(\alpha_{t})+ \delta \rho^{m_t}g_t^Td_t\\
\mid g_{t+1}^Td_t\mid\leq- \sigma g_t^Td_t
\end{array}\right..\eqno(2.4)$$
Set $\beta_t=\rho^{m_t},\ \gamma_t=\alpha_{t}(:),\ \gamma_{t+1}=\gamma_{t}+\beta_td_t,\ \alpha_{t+1}
=reshape(\gamma_{t+1}, n, k-1).$\\
\indent Step 5. Set $t:=t+1$. Go to step 2.\\
\vskip2mm
\hrule\vskip4mm

\noindent{\bf Remark 2.3.} To implement Algorithm 2.1, we first need to create three matlab files,
{\bf fun} file, {\bf gfun} file and {\bf frac} file, where the {\bf fun} file is used to compute
$F(\alpha _t)$, the {\bf gfun} file is used to calculate $\nabla F(\alpha_t)$, and the {\bf frac}
file is used to minimize $F(\alpha)$. In addition, the function $\alpha_{t}(:)$ returns the $n$
by $k-1$ vector $\gamma_t$ whose elements are taken column-wise from the matrix $\alpha_{t}$,
and the function $reshape(\gamma_{t+1}, n, k-1)$ returns the $n$ by $k-1$ matrix $\alpha_{t+1}$
whose elements are taken column-wise from $\gamma_{t+1}$.

By Theorem 4.3.5 [19, P.203], we can establish the global convergence theorem for Algorithm 2.1.

\noindent{\bf Theorem 2.3}. Suppose the function $F(\alpha)$ is twice continuous and differentiable,
the level set $$\Omega(\alpha_0)=\{\alpha\in R^{n\times (k-1)}\mid\ F(\alpha)\leq F(\alpha_0)\}$$ is
bounded, and the step length $\beta_t$ is generated by (2.4), where $\delta<\sigma<0.5$. Then the
sequence $\{\alpha_t\}$ generated by Algorithm 2.1 is guaranteed to globally converge, that is,
$$\lim\limits_{t\rightarrow \infty}inf\ \parallel\nabla F(\alpha_t)\parallel_F=0.$$

\vskip 3mm \noindent{\normalsize\bf 3. Numerical Experiments }
\vskip2mm

In this section, we use two numerical examples to illustrate that Algorithm 2.1 is feasible to
solve problem (2.1). All experiments are tested in \emph{Matlab R2010a}. We denote the relative residual error
$$\epsilon(t)=\frac {\sum\limits_{d=1}^{m}\|A^{(d)}-Y_t\|_F^2}{\sum\limits_{d=1}^{m} \|A^{(d)}\|_F^2},$$
and the gradient norm
$$\|g_t\|_F=\|\nabla F(\alpha_t)\|_F,$$
where $\alpha_t$ is the $t$th iterative matrix of Algorithm 2.1. We use the stopping criterion
$$\|g_t\|_F<1.0\times 10^{-4}.$$ And we choose the random matrix $rand(m,n)$ as the initial value in the
following examples, where the random matrix is generated by the Matlab function $rand(m,n).$

\noindent{\bf Example 3.1}. Consider problem (2.1) with $m=1$ and
{\scriptsize
$$ A=\left[\begin{array}{ccccc}
    1.0000  &  0.1849 &  -0.2867 &  -0.2997\\
    0.1849  &  1.0000  &  0.2851 &   0.2582\\
   -0.2867  &  0.2851  &  1.0000 &  -0.3100\\
   -0.2997  &  0.2582  & -0.3100 &   1.0000
\end{array}\right].$$}
\noindent{\bf Case I:} Set k=3. We use Algorithm 2.1 with the initial value
{\scriptsize
$$ \alpha_0=\left[\begin{array}{cc}
    0.0344  &  0.7952\\
    0.4387  &  0.1869\\
    0.3816  &  0.4898\\
    0.7655  &  0.4456
\end{array}\right]$$}
to solve problem (2.1). After 15 iterations, we get the solution $\widehat{\alpha}$ of problem (2.1)
{\scriptsize
$$\widehat{\alpha} \approx \alpha_{15}
=\left[\begin{array}{cc}
   -1.6439  &  1.3217\\
    1.2743  & -0.5270\\
    0.3266  &  0.8184\\
    2.2501  &  0.3043
\end{array}\right].$$ }
Hence, the solution $\widehat{Y}$ of problem (1.1) is
{\scriptsize
$$ \widehat{Y}=\left[\begin{array}{cccc}
    1.0000  &  0.2403  & -0.3495  & -0.3619\\
    0.2403  &  1.0000  &  0.3453  &  0.3179\\
   -0.3495  &  0.3453  &  1.0000  & -0.3777\\
   -0.3619  &  0.3179  & -0.3777  &  1.0000
\end{array}\right].$$}
And the curves of the relative residual error $\epsilon(t)$ and the gradient norm
$\|\nabla F(\alpha_t)\|_{F}$ are in Fig. 1.
\begin{figure}[htbp]
 \begin{minipage}{1\textwidth}
 \def\figurename{\footnotesize Fig.}
 \centering
\resizebox{14cm}{5.5cm}{\includegraphics{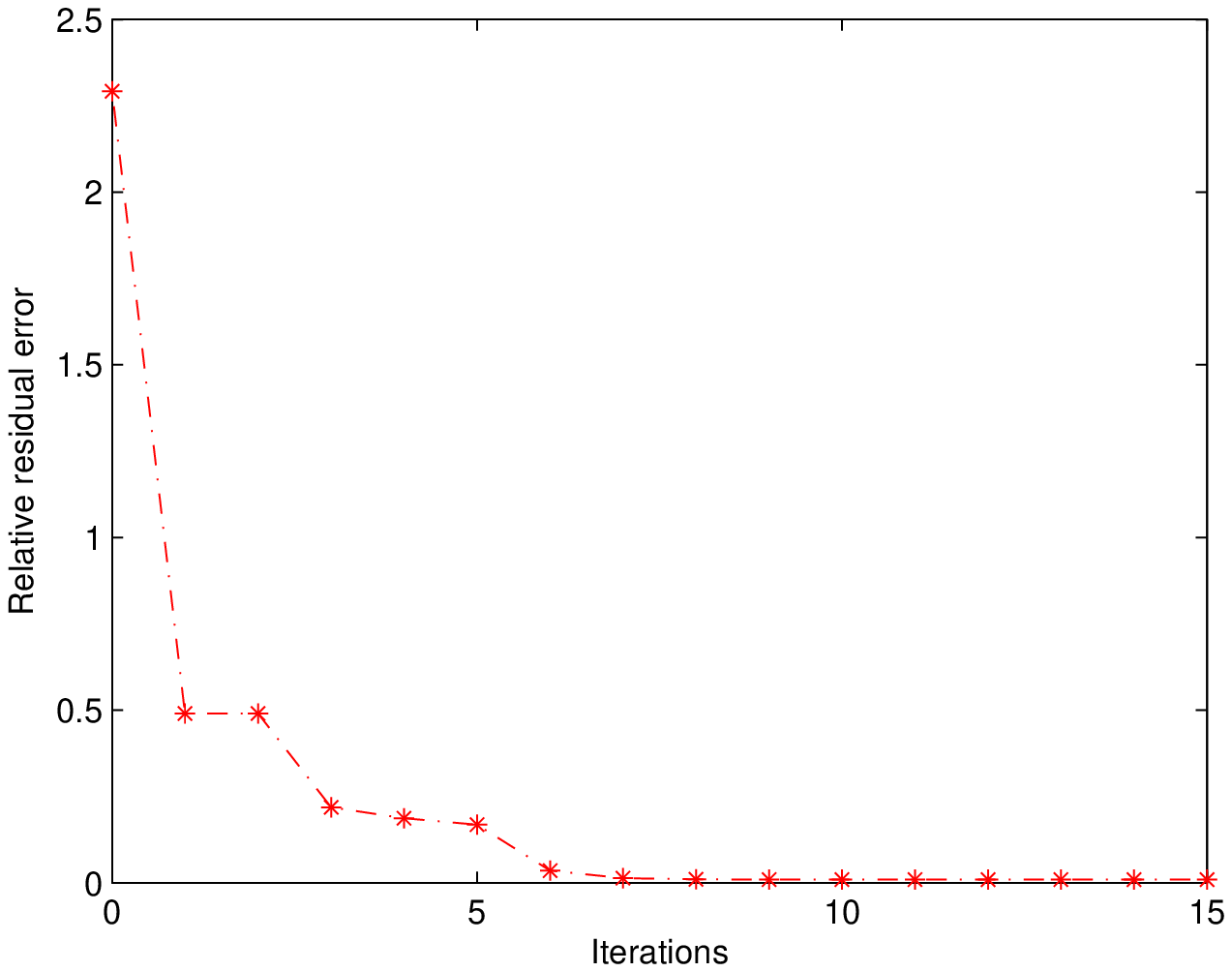}\includegraphics{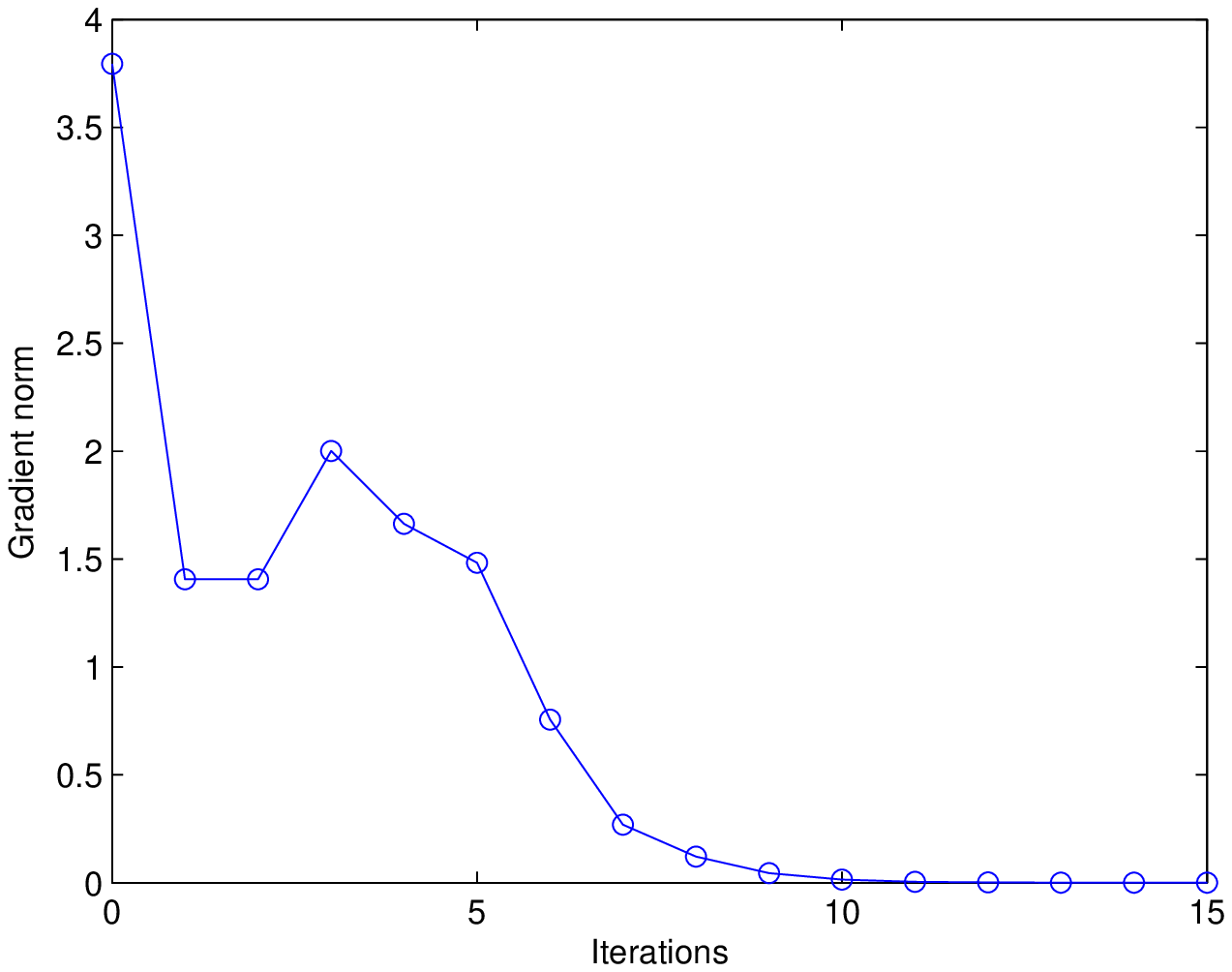}}
\caption{\footnotesize Convergence curves of the relative residual error
$\varepsilon(t)$ and the gradient norm $\|\nabla F(\alpha_t)\|_{F}$.}
   \end{minipage}
\end{figure}

\noindent{\bf Case II:} Set k=2. We use Algorithm 2.1 with the initial value
{\scriptsize
\begin{eqnarray*}
\alpha_0=\left[\begin{array}{c}
    0.9572\\
    0.4854\\
    0.8003\\
    0.1419
\end{array}\right]
\end{eqnarray*}}
to solve problem (2.1). After 13 iterations, we get the solution $\widehat{\alpha}$ of problem (2.1)
{\scriptsize
$$\widehat{\alpha}\approx\alpha_{13}=\left[\begin{array}{c}
    2.2975\\
    0.4993\\
    0.6773\\
   -1.0893
\end{array}\right].$$ }
Hence, the solution $\widehat{Y}$ of problem (1.1) is
{\scriptsize
$$ \widehat{Y}=\left[\begin{array}{cccc}
    1.0000  & -0.2254 &  -0.0494  & -0.9701\\
   -0.2254  &  1.0000  &  0.9842  & -0.0178\\
   -0.0494   & 0.9842  &  1.0000  & -0.1946\\
   -0.9701  & -0.0178  & -0.1946  &  1.0000
\end{array}\right].$$}
And the curves of the relative residual error $\epsilon(t)$ and the gradient norm
$\|\nabla F(\alpha_t)\|_{F}$ are in Fig. 2.
\begin{figure}[htbp]
 \begin{minipage}{1\textwidth}
 \def\figurename{\footnotesize Fig.}
 \centering
\resizebox{14cm}{5.5cm}{\includegraphics{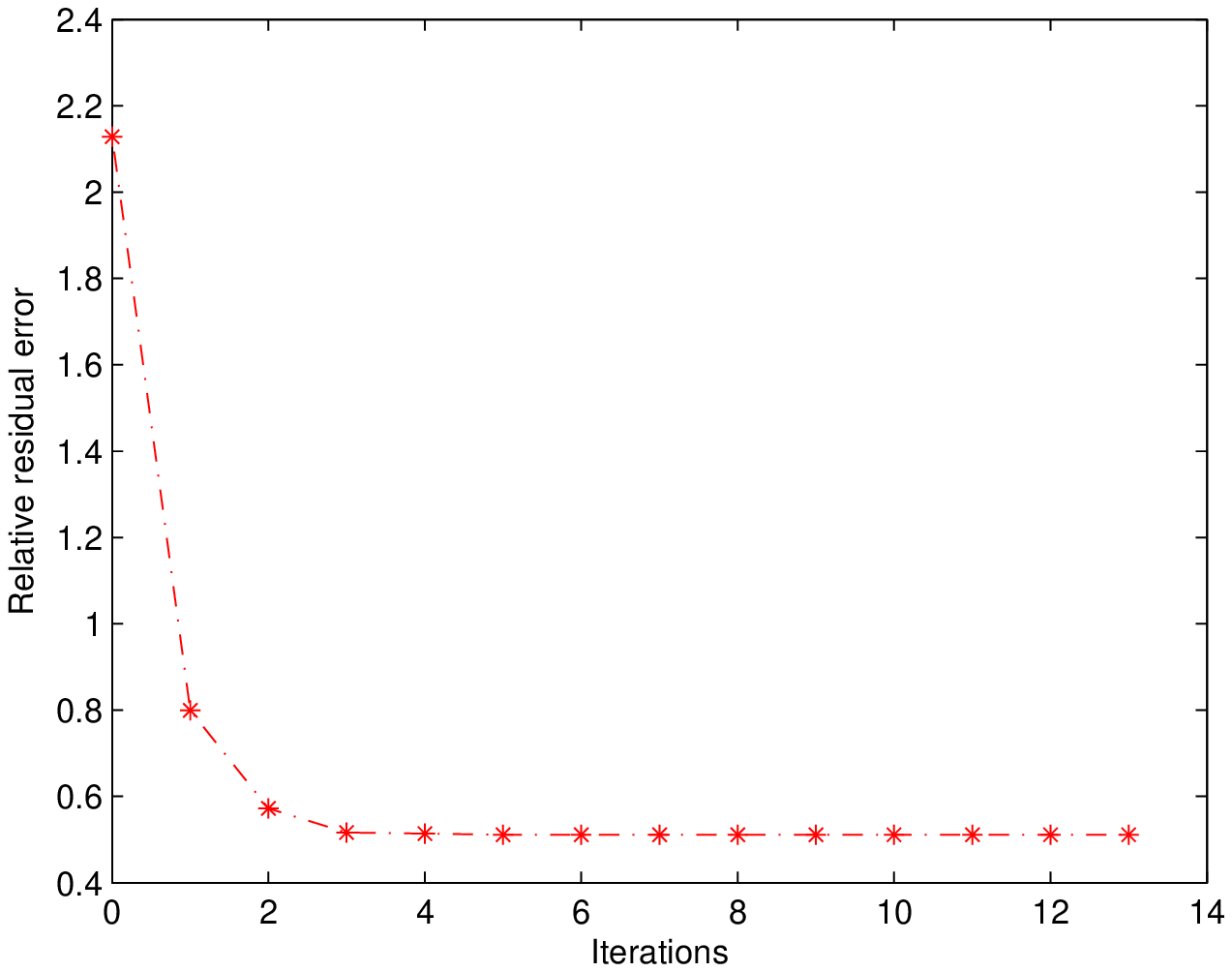}\includegraphics{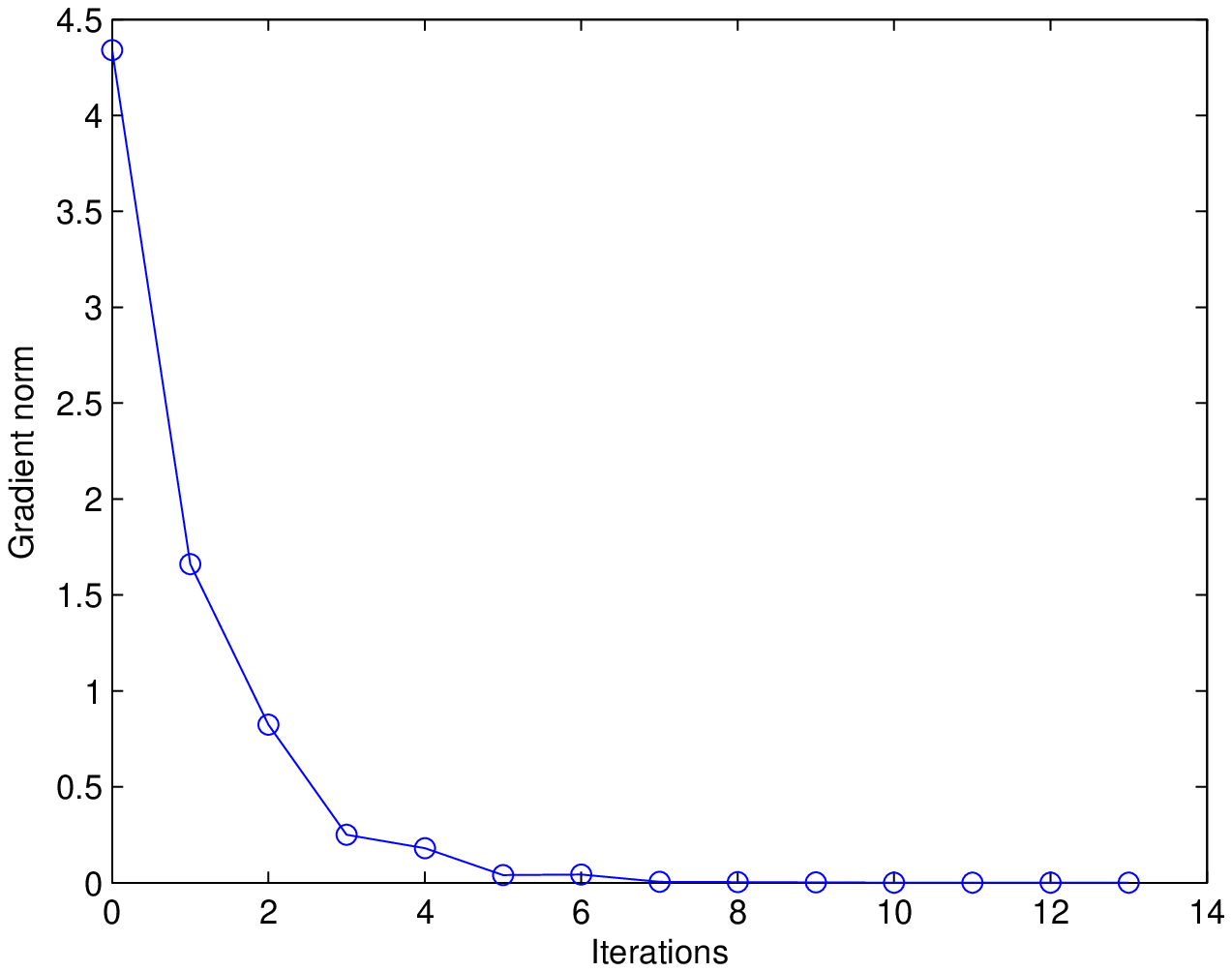}}
\caption{\footnotesize Convergence curves of the relative residual error
$\varepsilon(t)$ and the gradient norm $\|\nabla F(\alpha_t)\|_{F}$.}
   \end{minipage}
\end{figure}

In order to compare our algorithm with the Major algorithm in [11], we use them to solve
problem (2.1) with the same initial value. We list the number of iteration (denoted by "{\bf IT}"),
CPU time (denoted by "{\bf CPU}"), the gradient norm (denoted by "{\bf GN}") and the relative
residual error (denoted by "{\bf ERR}") in Table 1.
{\small
$$
\begin{tabular}{|c|c|c|c|c|}
\hline
\multicolumn{1}{|c|}{$rank\ \ k$} & \multicolumn{2}{|c|}{$2$}& \multicolumn{2}{|c|}{$3$} \\
 \hline
$Algorithm$ & $2.1$  & $Major$ & $2.1$ & $Major$  \\
\hline
$IT$        & $13$    & $14$ & $15$ & $19$  \\
\hline
$CPU (s)$   & $0.0312$  & $0.0624$ & $0.0468 $ & $0.1092$ \\
\hline
$GN$  & $7.2893\times 10^{-5}$  & $7.7516\times 10^{-5}$ & $7.2330\times 10^{-5}$ & $8.4253\times 10^{-5}$ \\
\hline
$ERR$ & $0.5111$ & $0.5168$ & $0.0092$ & $0.0105$ \\
\hline
\end{tabular}
$$
$$
Table  \ 1: Results \ for \ Example \ 3.1 \ with \ different \ values \ of \ rank  \ k
$$
}

Example 3.1 shows that Algorithm 2.1 is feasible to solve problem (1.1). Especially, Table 1
shows that our algorithm outperforms the Major algorithm [11] in both iterations and CPU time,
which indicates that our algorithm has faster convergence rate than the Major algorithm.

Nextly, we will use an example to show that our algorithm can be used to solve the generalized
low rank approximation of correlation matrices arising in the asset portfolio.

\noindent{\bf Example 3.2}. It is an important issue to calculate the more exact correlation matrix of assets in the portfolio selection.
For instance, suppose that an investor uses one unit money to buy a total of $11$ assets at the beginning of one period. There is a
relationship between any two assets of the portfolio because the price of each asset is related to some common factors in the financial market. The correlation matrix is one of the methods measuring the relation between assets. However, how to accurately compute the
correlation matrix is the key problem for the investor since the optimal investment policies is affected by the uncertainty of parameters in the correlation matrix. The daily price data of each asset in the portfolio are taken from the Wind database, which is a Chinese financial database, in order to obtain the correlation matrix. Five sets of the daily data are got by the sampling based on five different periods of the data.  Using the Matlab software, five correlation matrix of the eleven assets are given as follows.

{\tiny
$$A^{(1)}=\left[\begin{array}{ccccccccccc}
    1.0000 &  0.6712  &  0.5141  &  0.7085  &  0.9411  &  0.9435  &  0.9619 & 0.8106 &   0.5186 &  -0.0071  &  0.9514\\
    0.6712 &  1.0000  &  0.7421  &  0.7707  &  0.5058  &  0.5926  &  0.6942 & 0.7540 &   0.7738 &   0.5590  &  0.6122\\
    0.5141 &  0.7421  &  1.0000  &  0.4919  &  0.3912  &  0.3549  &  0.4227 & 0.4881 &   0.6179 &   0.4515  &  0.3700\\
    0.7085 &  0.7707  &  0.4919  &  1.0000  &  0.5708  &  0.7849  &  0.7084 & 0.6832 &   0.4142 &   0.1868  &  0.7442\\
    0.9411 &  0.5058  &  0.3912  &  0.5708  &  1.0000  &  0.8967  &  0.9175 & 0.6512 &   0.3372 &  -0.2023  &  0.9251\\
    0.9435 &  0.5926  &  0.3549  &  0.7849  &  0.8967  &  1.0000  &  0.9316 & 0.7522 &   0.3542 &  -0.1386  &  0.9618\\
    0.9619 &  0.6942  &  0.4227  &  0.7084  &  0.9175  &  0.9316  &  1.0000 & 0.8441 &   0.5710 &   0.0352  &  0.9483\\
    0.8106 &  0.7540  &  0.4881  &  0.6832  &  0.6512  &  0.7522  &  0.8441 & 1.0000 &   0.8176 &   0.3378  &  0.7849\\
    0.5186 &  0.7738  &  0.6179  &  0.4142  &  0.3372  &  0.3542  &  0.5710 & 0.8176 &   1.0000 &   0.6533  &  0.4024\\
   -0.0071 &  0.5590  &  0.4515  &  0.1868  & -0.2023  & -0.1386  &  0.0352 & 0.3378 &   0.6533 &   1.0000  & -0.0495\\
    0.9514 &  0.6122  &  0.3700  &  0.7442  &  0.9251  &  0.9618  &  0.9483 & 0.7849 &   0.4024 &  -0.0495  &  1.0000
\end{array}\right],
$$
$$A^{(2)}=\left[\begin{array}{ccccccccccc}
    1.0000   & 0.8140  &  0.9019  &  0.8838   & 0.4088  &  0.9100  &  0.2976 &  0.5686  &  0.2685  &  0.6239  &  0.3775\\
    0.8140   & 1.0000  &  0.9564  &  0.9502   & 0.5080  &  0.8474  &  0.3311 &  0.3952  &  0.1728  &  0.8082  &  0.4650\\
    0.9019   & 0.9564  &  1.0000  &  0.9586   & 0.5407  &  0.9202  &  0.2926 &  0.4808  &  0.1827  &  0.8114  &  0.5034\\
    0.8838   & 0.9502  &  0.9586  &  1.0000   & 0.4479  &  0.8992  &  0.4163 &  0.5251  &  0.2953  &  0.7441  &  0.4152\\
    0.4088   & 0.5080  &  0.5407  &  0.4479   & 1.0000  &  0.4441  & -0.3338 & -0.1371  & -0.3909  &  0.7998  &  0.9451\\
    0.9100   & 0.8474  &  0.9202  &  0.8992   & 0.4441  &  1.0000  &  0.3586 &  0.5868  &  0.2569  &  0.7422  &  0.4388\\
    0.2976   & 0.3311  &  0.2926  &  0.4163   &-0.3338  &  0.3586  &  1.0000 &  0.8315  &  0.9211  &  0.0386  & -0.3588\\
    0.5686   & 0.3952  &  0.4808  &  0.5251   &-0.1371  &  0.5868  &  0.8315 &  1.0000  &  0.8911  &  0.2210  & -0.1591\\
    0.2685   & 0.1728  &  0.1827  &  0.2953   &-0.3909  &  0.2569  &  0.9211 &  0.8911  &  1.0000  & -0.0558  & -0.4198\\
    0.6239   & 0.8082  &  0.8114  &  0.7441   & 0.7998  &  0.7422  &  0.0386 &  0.2210  & -0.0558  &  1.0000  &  0.7942\\
    0.3775   & 0.4650  &  0.5034  &  0.4152   & 0.9451  &  0.4388  & -0.3588 & -0.1591  & -0.4198  &  0.7942  &  1.0000
\end{array}\right],
$$
$$A^{(3)}=\left[\begin{array}{ccccccccccc}
    1.0000  &  0.8581 &   0.8033 &   0.7763  &  0.5692 &   0.8994  & -0.0383 & -0.1388  & -0.2484  &  0.7421  &  0.5445\\
    0.8581  &  1.0000 &   0.8446 &   0.7744  &  0.4408 &   0.8166  &  0.1116 & -0.1725  & -0.1207  &  0.5586  &  0.3944\\
    0.8033  &  0.8446 &   1.0000 &   0.8788  &  0.2731 &   0.8565  &  0.2448 & -0.0567  &  0.1683  &  0.4772  &  0.2438\\
    0.7763  &  0.7744 &   0.8788 &   1.0000  &  0.3428 &   0.8868  &  0.2869 &  0.0620  &  0.2111  &  0.4601  &  0.3225\\
    0.5692  &  0.4408 &   0.2731 &   0.3428  &  1.0000 &   0.4730  & -0.5636 & -0.4667  & -0.6824  &  0.8637  &  0.9721\\
    0.8994  &  0.8166 &   0.8565 &   0.8868  &  0.4730 &   1.0000  &  0.1251 & -0.0813  & -0.0267  &  0.6438  &  0.4551\\
   -0.0383  &  0.1116 &   0.2448 &   0.2869  & -0.5636 &   0.1251  &  1.0000 &  0.6858  &  0.8411  & -0.5392  & -0.5661\\
   -0.1388  & -0.1725 &  -0.0567 &   0.0620  & -0.4667 &  -0.0813  &  0.6858 &  1.0000  &  0.7263  & -0.4975  & -0.4254\\
   -0.2484  & -0.1207 &   0.1683 &   0.2111  & -0.6824 &  -0.0267  &  0.8411 &  0.7263  &  1.0000  & -0.6348  & -0.6618\\
    0.7421  &  0.5586 &   0.4772 &   0.4601  &  0.8637 &   0.6438  & -0.5392 & -0.4975  & -0.6348  &  1.0000  &  0.8715\\
    0.5445  &  0.3944 &   0.2438 &   0.3225  &  0.9721 &   0.4551  & -0.5661 & -0.4254  & -0.6618  &  0.8715  &  1.0000
\end{array}\right],
$$
$$A^{(4)}=\left[\begin{array}{ccccccccccc}
    1.0000 & 0.6803 & 0.7064 & 0.8565 &-0.2759 & 0.5470 & 0.4280 & 0.3874 & 0.3382 & 0.3684 &-0.2266\\
    0.6803 & 1.0000 & 0.7341 & 0.7650 &-0.2123 & 0.7590 &-0.1643 &-0.1412 &-0.1483 &-0.0227 &-0.1681\\
    0.7064 & 0.7341 & 1.0000 & 0.7334 &-0.2411 & 0.5976 &-0.0299 &-0.0849 &-0.1307 & 0.0605 &-0.1856\\
    0.8565 & 0.7650 & 0.7334 & 1.0000 &-0.2705 & 0.6115 & 0.2210 & 0.1977 & 0.1355 & 0.2755 &-0.1968\\
   -0.2759 &-0.2123 &-0.2411 &-0.2705 & 1.0000 &-0.1890 &-0.1144 &-0.0014 & 0.0969 & 0.4612 & 0.9336\\
    0.5470 & 0.7590 & 0.5976 & 0.6115 &-0.1890 & 1.0000 &-0.3366 &-0.2152 &-0.2045 &-0.2603 &-0.1309\\
    0.4280 &-0.1643 &-0.0299 & 0.2210 &-0.1144 &-0.3366 & 1.0000 & 0.8938 & 0.8434 & 0.6356 &-0.1117\\
    0.3874 &-0.1412 &-0.0849 & 0.1977 &-0.0014 &-0.2152 & 0.8938 & 1.0000 & 0.9486 & 0.6122 & 0.0158\\
    0.3382 &-0.1483 &-0.1307 & 0.1355 & 0.0969 &-0.2045 & 0.8434 & 0.9486 & 1.0000 & 0.5966 & 0.1128\\
    0.3684 &-0.0227 & 0.0605 & 0.2755 & 0.4612 &-0.2603 & 0.6356 & 0.6122 & 0.5966 & 1.0000 & 0.5056\\
   -0.2266 &-0.1681 &-0.1856 &-0.1968 & 0.9336 &-0.1309 &-0.1117 & 0.0158 & 0.1128 & 0.5056 & 1.0000
\end{array}\right].
$$
$$
A^{(5)}=\left[\begin{array}{ccccccccccc}
    1.0000 & 0.2118 & 0.1238 & 0.2178 &-0.2533 &-0.0778 & 0.7000 & 0.3288 & 0.1310 &-0.0052 & 0.1428\\
    0.2118 & 1.0000 & 0.8882 & 0.7828 & 0.6747 &-0.8135 & 0.3794 & 0.8962 & 0.8687 & 0.6974 & 0.4794\\
    0.1238 & 0.8882 & 1.0000 & 0.6828 & 0.7155 &-0.9202 & 0.4205 & 0.7974 & 0.9306 & 0.8604 & 0.7235\\
    0.2178 & 0.7828 & 0.6828 & 1.0000 & 0.6836 &-0.5435 & 0.3370 & 0.6787 & 0.6683 & 0.3548 & 0.1678\\
   -0.2533 & 0.6747 & 0.7155 & 0.6836 & 1.0000 &-0.6628 & 0.0448 & 0.4736 & 0.6978 & 0.5897 & 0.3092\\
   -0.0778 &-0.8135 &-0.9202 &-0.5435 &-0.6628 & 1.0000 &-0.4037 &-0.7538 &-0.8888 &-0.8936 &-0.7417\\
    0.7000 & 0.3794 & 0.4205 & 0.3370 & 0.0448 &-0.4037 & 1.0000 & 0.5818 & 0.4775 & 0.3655 & 0.4722\\
    0.3288 & 0.8962 & 0.7974 & 0.6787 & 0.4736 &-0.7538 & 0.5818 & 1.0000 & 0.8544 & 0.6521 & 0.5163\\
    0.1310 & 0.8687 & 0.9306 & 0.6683 & 0.6978 &-0.8888 & 0.4775 & 0.8544 & 1.0000 & 0.8203 & 0.6500\\
   -0.0052 & 0.6974 & 0.8604 & 0.3548 & 0.5897 &-0.8936 & 0.3655 & 0.6521 & 0.8203 & 1.0000 & 0.8810\\
    0.1428 & 0.4794 & 0.7235 & 0.1678 & 0.3092 &-0.7417 & 0.4722 & 0.5163 & 0.6500 & 0.8810 & 1.0000
\end{array}\right],
$$}

Set k=3, and we use Algorithm 2.1 with the initial value
{\scriptsize
$$\alpha_0=\left[\begin{array}{cc}
    0.0462  &  0.1869\\
    0.0971  &  0.4898\\
    0.8235  &  0.4456\\
    0.6948  &  0.6463\\
    0.3171  &  0.7094\\
    0.9502  &  0.7547\\
    0.0344  &  0.2760\\
    0.4387  &  0.6797\\
    0.3816  &  0.6551\\
    0.7655  &  0.1626\\
    0.7952  &  0.1190
\end{array}\right]
$$}to solve problem (2.1). After 57 iterations, we get the solution $\widehat{\alpha}$ of problem (2.1)
{\scriptsize
$$\widehat{\alpha}\approx \alpha_{57}=\left[\begin{array}{cc}
    0.4179  &  1.2147\\
    0.4126  &  0.4239\\
    0.3730   & 0.3196\\
    0.2868  &  0.9097\\
    1.2956  & -0.3683\\
    0.9810  &  1.3296\\
   -0.7615  &  0.6138\\
   -0.7709  &  0.9219\\
   -0.8959  &  0.9181\\
    1.0963  & -0.7234\\
    1.2672  & -0.4347
\end{array}\right].
$$}
Hence, the solution $\widehat{Y}$ of problem (1.1) is
{\scriptsize
$$ \widehat{Y}=\left[\begin{array}{ccccccccccc}
    1.0000  &  0.9517  &  0.9436  &  0.9861 &   0.2436 &   0.8434  &  0.4306  &  0.3849  &  0.2680  &  0.2879  &  0.2428\\
    0.9517  &  1.0000  &  0.9984  &  0.9790 &   0.5200  &  0.7152  &  0.3913  &  0.4117  &  0.2967  &  0.5651  &  0.5239\\
    0.9436  &  0.9984  &  1.0000  &  0.9789  &  0.5240  &  0.6790  &  0.4334  &  0.4588  &  0.3468  &  0.5887  &  0.5318\\
    0.9861  &  0.9790  &  0.9789  &  1.0000  &  0.3393  &  0.7482  &  0.5075  &  0.4909 &   0.3784  &  0.4225  &  0.3473\\
    0.2436  &  0.5200  &  0.5240  &  0.3393  &  1.0000  &  0.0498  & -0.1720  &  0.0093 &  -0.0410  &  0.9268  &  0.9976\\
    0.8434  &  0.7152  &  0.6790  &  0.7482  &  0.0498  &  1.0000  & -0.0301  & -0.1326 &  -0.2472  & -0.0888 &   0.0138\\
    0.4306  &  0.3913  &  0.4334  &  0.5075  & -0.1720  & -0.0301  &  1.0000  &  0.9773  &  0.9662  &  0.1886 &  -0.1121\\
    0.3849  &  0.4117  &  0.4588  &  0.4909  &  0.0093 &  -0.1326  &  0.9773  &  1.0000  &  0.9922 &   0.3739 &   0.0731\\
    0.2680   & 0.2967  &  0.3468  &  0.3784  & -0.0410 &  -0.2472  &  0.9662  &  0.9922  &  1.0000 &   0.3345 &   0.0256\\
    0.2879  &  0.5651  &  0.5887  &  0.4225  &  0.9268  & -0.0888  &  0.1886  &  0.3739  &  0.3345 &   1.0000 &   0.9503\\
    0.2428  &  0.5239  &  0.5318  &  0.3473  &  0.9976 &   0.0138  & -0.1121   & 0.0731  &  0.0256  &  0.9503 &   1.0000
\end{array}\right].
$$}
And the curves of the relative residual error $\epsilon(t)$ and the gradient norm $\|\nabla F(\alpha_t)\|_{F}$ are in Fig. 3.
\begin{figure}[htbp]
 \begin{minipage}{1\textwidth}
 \def\figurename{\footnotesize Fig.}
 \centering
\resizebox{14cm}{5.5cm}{\includegraphics{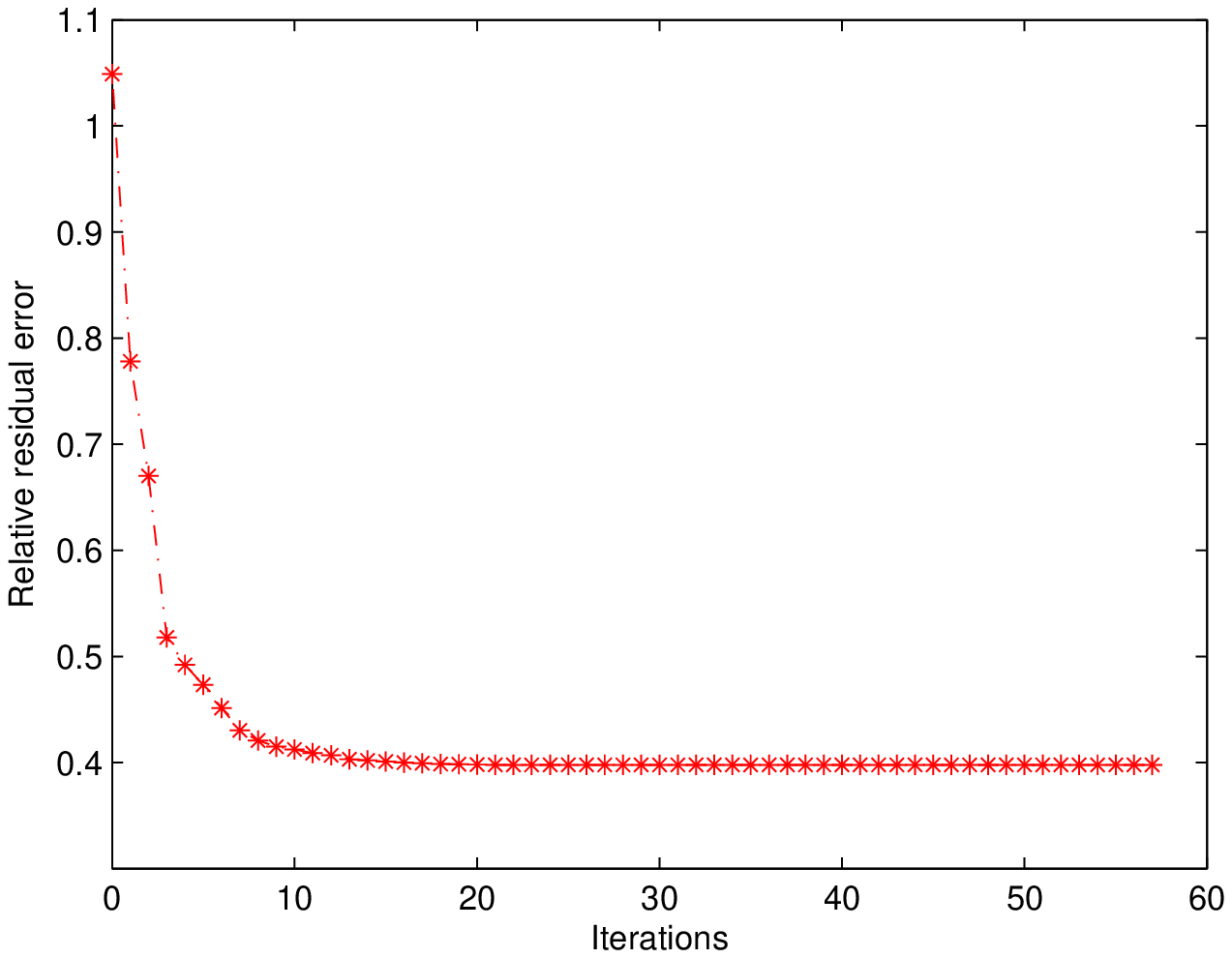}\includegraphics{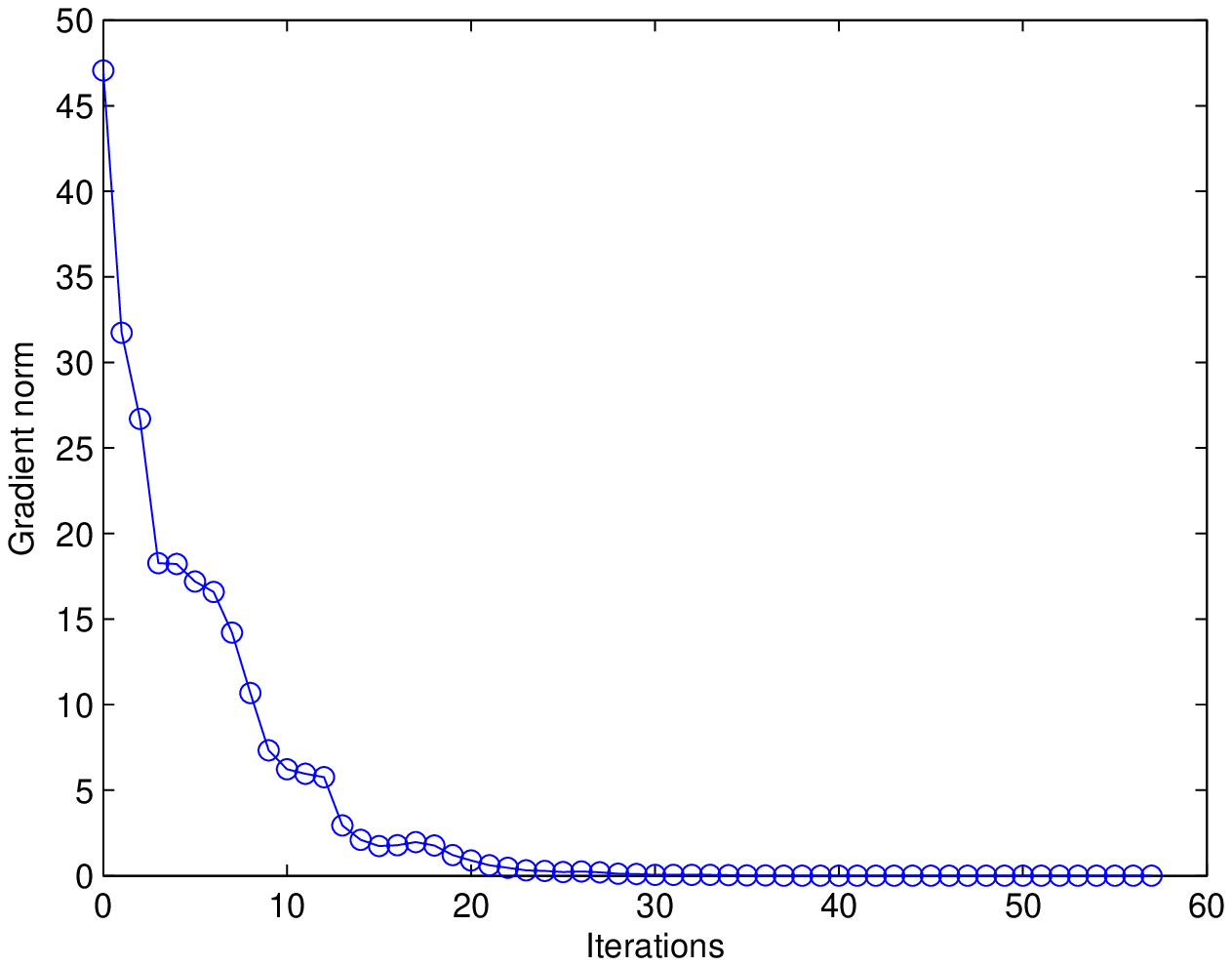}}
\caption{\footnotesize Convergence curves of the relative residual error
$\varepsilon(t)$ and the gradient norm $\|\nabla F(\alpha_t)\|_{F}$.}
   \end{minipage}
\end{figure}

For the above example, we use Algorithm 2.1 to solve problem (2.1) with different rank. We
list the number of iteration (denoted by "{\bf IT}") , CPU time (denoted by "{\bf CPU}"),
the gradient norm (denoted by "{\bf GN}") and the relative residual error (denoted by "{\bf ERR}") in Table 2.
{\small
$$
\begin{tabular}{|c|c|c|c|c|}
\hline
$rank\ k$ & $2$ & $3$  & $4$ & $5$  \\
\hline
$IT$ & $44$  & $57$ & $1005$ & $2000$ \\
\hline
$CPU (s)$ & $ 0.1404$  & $0.2184$ & $8.1121 $ & $21.9649$ \\
\hline
$GN$ & $5.2915\times 10^{-5}$  & $9.5734\times 10^{-5}$ & $9.9882\times 10^{-5}$ & $0.3687$ \\
\hline
$ERR$ & $0.5879$ & $0.3977$ & $0.4532$ & $0.4087$ \\
\hline
\end{tabular}
$$
$$Table\ 2: Results \ for  \ Example \ 3.2 \ with \ different \ rank   \ by \ Algorithm \ 2.1$$
}

Fig. 3 and Table 2 show that Algorithm 2.1 can be used to solve the generalized low rank
approximation of correlation matrices arising in the asset portfolio. What is more important,
when the investor uses the matrix $\widehat{Y}$ obtained by using Algorithm 2.1 to analyze the relationship
between any two assets, some noise in the data can be reduced because the correlation matrix of
assets is an important factor for selecting assets in portfolio.

\vskip 3mm \noindent{\normalsize\bf 4. Conclusion }
\vskip2mm
The generalized low rank approximation of correlation matrices is widely used in the asset
portfolio and risk management. It is a difficult matrix optimization problem, and the
difficulties lie in how to deal with its feasible set and complex structure. In this paper,
we use the Gramian representation together with special trigonometric function transform to
overcome these difficulties, and develop a new algorithm to solve it. Numerical examples
show that our new method  is feasible and effective. Moreover, the theory and algorithm of
this paper can be extended to solve the low rank approximation in Li-Qi [9], that is, the
nearest low rank approximation of a correlation matrix to the given symmetric matrix.

\vskip 3mm \noindent{\large\bf Acknowledgements}\vskip 2mm

The authors wish to thank Prof. Richard A Brualdi and the anonymous
referee for providing very useful suggestions for improving this
paper. The authors also thank Prof. Qingwen Wang for discussing
the properties of the objective function.

\end{document}